\documentclass[12pt,reqno]{amsart}
\usepackage{amsmath,amssymb,amscd,latexsym,eufrak, mathrsfs}
\usepackage{graphicx}
\newcommand{\svskip}{\vspace{3mm}}
\newcommand{\C}{{\mathbb C}}
\newcommand{\F}{{\mathbb F}}
\newcommand{\Q}{{\mathbb Q}}
\newcommand{\Bl}{{\rm Bl}}
\newcommand{\BP}{{\mathbb P}}
\newcommand{\gm}{\EuFrak{m}}

\newcommand{\QED}{{\unskip\nobreak\hfil\penalty50\quad\null\nobreak\hfil
{$\Box$}\parfillskip0pt\finalhyphendemerits0\par\medskip}}
\newcommand{\Proof}{\noindent{\bf Proof.}\quad}
\newcommand{\A}{{\mathbb A}}
\newcommand{\Spec}{{\rm Spec}\:}
\newcommand{\Proj}{{\rm Proj}}

\newcommand{\Aut}{{\rm Aut}}
\newcommand{\Pic}{{\rm Pic}\:}
\newcommand{\codim}{{\rm codim}\:}
\newcommand{\Supp}{{\rm Supp}}
\newcommand{\Ker}{{\rm Ker}\:}
\newcommand{\id}{{\rm id}}
\newcommand{\dto}{\dashrightarrow}
\newcommand{\SO}{{\mathcal O}}
\newcommand{\wt}{\widetilde}
\newcommand{\wh}{\widehat}
\newcommand{\ol}{\overline}

\newcommand{\GL}{{\rm GL}}
\newcommand{\PGL}{{\rm PGL}}
\newtheorem{thm}{Theorem}[section]
\newtheorem{lem}[thm]{Lemma}

\newtheorem{remark}[thm]{Remark}
\newtheorem{example}[thm]{Example}

\newtheorem{problem}[thm]{Problem}

\newcommand{\st}[1]{\stackrel{{#1}}{\longrightarrow}}

\begin{document}
\title{Equivariant Jacobian Conjecture in dimension two}

\author{Masayoshi Miyanishi}

\thanks{The author was supported by Grant-in-Aid for Scientific Research (C), No. 16K05115, JSPS}

\address{Research Center for Mathematics and Data Science, Kwansei Gakuin University \\
2-1 Gakuen, Sanda 669-1337, Japan}
\email{miyanisi@kwansei.ac.jp}

\maketitle

\begin{abstract}
Let $G$ be a small finite subgroup of $\GL(2,\C)$ and let $\wt{\varphi} : \A^2 \to \A^2$ be a $G$-equivariant 
\'etale endomorphism of the affine plane. We show that $\wt{\varphi}$ is an automorphism if the order of $G$ is even. 
The proof depends on an analysis of a quasi-\'etale endomorphism $\varphi$ induced by $\wt{\varphi}$ on the 
singular quotient surface $\A^2/G$ whose smooth part $X^\circ$ has the standard $\A^1_\ast$-fibration $p^\circ : X^\circ 
\to \BP^1$. If $\varphi$ preserves the standard $\A^1_\ast$-fibration $p^\circ$ then both $\varphi$ and $\wt{\varphi}$ are 
automorphisms. We look for the condition with which $\varphi$ preserves the standard $\A^1_\ast$-fibration and prove 
as a consequence that $\wt{\varphi}$ is an automorphism if $G$ has even order. 
\end{abstract}

\section*{Introduction}
Let $X^\circ$ be a smooth algebraic variety defined over the complex number field $\C$ and let $\varphi^\circ : 
X^\circ \to X^\circ$ be an \'etale endomorphism. The generalized Jacobian conjecture (GJC, for short) asks if 
$\varphi^\circ$ is a finite morphism. We considered this problem extensively in the case $\dim X^\circ=2$. 
See the references \cite{GM, MM, M1, KM}. The obtained results support the conjecture affirmatively, but 
there are scattered counterexamples sometimes with non-trivial families of non-finite \'etale endomorphisms \cite[Thm. C]{DP}. 
In \cite{GM}, we treated the case where $X^\circ$ is 
the smooth part of a singular normal affine surface $X$. Then the endomorphism $\varphi^\circ$ extends to 
an endomorphism $\varphi : X \to X$. An endomorphism $\varphi : X \to X$ is called a {\em quasi-\'etale endomorphism} 
if it induces an \'etale endomorphism $\varphi^\circ : X^\circ \to X^\circ$. In the subsequent arguments, we often 
denote $\varphi : X \to X$ by $\varphi : X_u \to X_l$ to distinguish the source $X$ from the target $X$, 
where ``$u$'' (resp. ``$l$'') means ``upper'' (resp. `` lower'').

In the present article, we consider the case $X=\A^2/G$, where $G$ is a nontrivial small finite subgroup of 
$\GL(2,\C)$. Hence $X$ has a unique quotient singularity $P_0$, which is the image of the origin $O$ in $\A^2$. 
A quasi-\'etale endomorphism $\varphi$ sends the point $P_0$ to itself. Let $\wh{X}$ be a minimal resolution of 
singularity of $X$. If $G$ is non-cyclic, then the exceptional locus $\Gamma$ is an admissible rational fork 
which consists of the central component $S_0$ and three linear chains $T_i~(i=1,2,3)$ sprouting out of $S_0$ such that 
the determinants $\delta_i=d(T_i)$ form a Platonic triplet. If $G$ is cyclic, then the exceptional locus $\Gamma$ is 
an admissible rational rod. For the terminology and the definitions in the theory of open algebraic surfaces, 
a general reference is \cite[p. 93]{M3}.

We first show that $\varphi$ lifts to an \'etale endomorphism $\wt{\varphi} : \A^2 \to \A^2$ such that 
$\wt{\varphi}^{-1}(O)=(O)$ and $\wt{\varphi}(gx)=\chi(g)\wt{\varphi}(x)$ for $g \in G$ and $x \in \A^2$, where 
$\chi$ is a group automorphism of $G$. We call $\wt{\varphi}$ a $(G,\chi)$-equivariant \'etale endomorphism. 
Since $G$ is a finite group, there exists an integer $N$ such that $\chi^N=\id$. Then $\varphi^N$ is a 
$G$-equivariant endomorphism, that is, $\wt{\varphi}^N(gx)=g\wt{\varphi}^N(x)$. 

We show that the quasi-\'etale endomorphism $\varphi$ extends to an \'etale endomorphism $\wh{\varphi} : 
\wh{X} \to \wh{X}$ provided $\varphi$ lifts up to a $G$-equivariant endomorphism $\wt{\varphi}$. Under the same 
condition, $\varphi$ is also \'etale at the singular point $P_0$ (see Lemma \ref{Lemma 2.2}). 

We make here the following two well-known remarks for the case of an \'etale endomorphism $\wt{\varphi} : \A^2 \to \A^2$ 
which is equivariant with respect to a (not necessarily small) finite subgroup $G$ of $\GL(2,\C)$.

(1)~Let $N$ be the normal subgroup of $G$ consisting of all peseudo-reflections of $G$. Then the algebraic quotient $\A^2/N$ 
is isomorphic to the affine plane $\A^2$, and $\wt{\varphi}$ induces a $G/N$-equivariant \'etale endomorphism 
$\wt{\varphi}/N : \A^2 \to \A^2$, where $G/N$ is a small finite subgroup of $\GL(2,\C)$.

(2)~Let $X=(\A^2/N)/(G/N)$. Then $\wt{\varphi}$ induces a quasi-\'etale endomorphism $\varphi$ of $X$, and $\wt{\varphi}$ 
is an automorphism if and only if so is $\varphi$. In fact, if $\varphi$ is an automorphism, then $\wt{\varphi}$ is a 
finite \'etale covering which is an automorphism because $\A^2$ is simply-connected.

By the second remark, we can restrict ourselves to the case where $G$ is a small finite subgroup.

In section one, we define the standard $\A^1_\ast$-fibration $p^\circ : X^\circ \to C$ on the smooth part
$X^\circ=X\setminus\{P_0\}$, where $C \cong \BP^1$, and give the logarithmic canonical divisor formula in terms 
of the standard fibration. 

In section two, we prove that if a quasi-\'etale endomorphism $\varphi$ of $X$ preserves the standard 
$\A^1_\ast$-fibration $p^\circ$ then $\varphi$ is an automorphism (see Theorem \ref{Theorem 2.4}). 

In section three, we consider the conditions under which the given \'etale endomorphism $\varphi$ 
of $X$ preserves the $\A^1_\ast$-fibration $p^\circ$ of $X^\circ$ (see Lemmas \ref{Lemma 3.5} and \ref{Lemma 3.7}).

In section four, we show that any quasi-\'etale endomorphism $\varphi$ preserves the standard $\A^1_*$-fibration if 
$G$ is cyclic and the action is given by ${}^\zeta(x,y)=(\zeta x, \zeta y)$ if $G=\langle\zeta\rangle$. Hence 
$\varphi$ is an automorphism in this case.

The last observation gives a surprising result. This remark is due to R.V. Gurjar to whom the author expresses his 
deep gratitude. 
\svskip

\noindent
{\bf Theorem.}~~{\em Let $G$ be a small finite subgroup of $\GL(2,\C)$ and let $X=\A^2/G$. Let $\varphi : X \to X$ 
be a quasi-\'etale endomorphism. If $G$ has even order then $\varphi$ is an automorphism. In particular, if 
$\wt{\varphi} : \A^2 \to \A^2$ is a $G$-equivariant \'etale endomorphism then $\wt{\varphi}$ is an automorphism.}
\svskip

As an example, the above theorem implies that a polynomial endomorphism $(x,y) \mapsto (f(x,y), g(x,y))$ of the affine 
plane $\A^2$ is an automorphism if the Jacobian determinant of $f, g$ is a nonzero constant and $f(-x,-y)=-f(x,y)$ and 
$g(-x,-y)=-g(x,y)$. The case where $G$ is cyclic and $d > 1$ (see Section one for $d$) is regrettably
 unsettled. 

We use the terms like {\em subdivisional} or {\em sprouting} blowing-up of the boundary divisor for which the readers are 
referred to \cite[p. 87]{M3}. To accord with the current use, a smooth normal completion (see \cite[p. 65]{M3} is called 
a {\em log smooth completion}.

The author thanks the referees for persevering critical reading of the manuscript and many suggestions of improvement.

\section{The standard $\A^1_\ast$-fibration on $\A^2/G$}

Let $X=\A^2/G$ and let $X^\circ=X\setminus\{P_0\}$, where $P_0$ is the unique singular point and the image of the origin 
$O$ in $\A^2$ by the quotient morphism $\pi : \A^2 \to \A^2/G$. Suppose first that $G$ is non-cyclic. Then $X^\circ$ is embedded 
as an open set into a smooth projective $\BP^1$-fiber space $p : V \to C \cong \BP^1$ as $X^\circ=V-(S_0+T_1+T_2+T_3)-(S_1+R_1+R_2+R_3)$, 
where :
\begin{enumerate}
\item[(1)]
$S_0$ and $S_1$ are sections of the $\BP^1$-fibration $p$.
\item[(2)]
$S_0+T_1+T_2+T_3$ is the exceptional locus of the minimal resolution of the singular point $P_0$. Hence it is an 
admissible rational fork with negative definite intersection matrix.
\item[(3)]
There are three $(-1)$-curves $F_1, F_2, F_3$ such that $T_i+F_i+R_i$ is a linear chain and supports a degenerate $\BP^1$-fiber 
$\ell_i$ of $p$ with the unique $(-1)$-curve $F_i$ for $i=1, 2, 3$. Let $m_i$ be the multiplicity of $F_i$ in $\ell_i$. 
Then the set $\{m_1,m_2,m_3\}$ is one of the Platonic triplets $\{2,2,n\} (n \ge 2), \{2,3,3\}, \{2,3,4\}$ and $\{2,3,5\}$. 
The linear part $R_i$ in $\ell_i$ is uniquely determined by $T_i+F_i$.
\item[(4)]
There are no singular fibers of $p$ other than the above three fibers $\ell_i$.
\end{enumerate}

The restriction $p^\circ=p\vert_{X^\circ} : X^\circ \to C$ defines a Platonic $\A^1_\ast$-fiber space. For the details, see 
\cite[p.143]{M2}.
\svskip

Suppose next that $G$ is a cyclic group of order $n$, which we identify with the group $G=\{\zeta^i \mid 0 \le i < n\}$ with 
a primitive $n$th root $\zeta$ of $1$. We choose a system of coordinates $\{x,y\}$ on $\A^2$ so that $G$ acts by 
${}^\zeta(x,y)=(\zeta x, \zeta^d y)$, where $0 < d < n$ and $\gcd(n,d)=1$. We say that the $G$-action has type $(n,d)$. 
Then the minimal resolution of singularity at $P_0$ has the exceptional locus indicated by the dual graph

\raisebox{-15mm}{
\begin{picture}(110,25)(-20,0)
\setlength{\unitlength}{0.9mm}
\put(5,10){\circle{2}}
\put(6,10){\line(1,0){18}}
\put(25,10){\circle{2}}
\put(26,10){\line(1,0){18}}
\put(45,10){\circle{2}}
\put(46,10){\line(1,0){19}}
\put(70,10){$\dots\dots$}
\put(85,10){\line(1,0){19}}
\put(105,10){\circle{2}}
\put(5,5){$E_1$}
\put(4,15){$-a_1$}
\put(25,5){$E_2$}
\put(24,15){$-a_2$}
\put(45,5){$E_3$}
\put(44,15){$-a_3$}
\put(105,5){$E_r$}
\put(104,15){$-a_r$}
\end{picture}}

\noindent
where $E_i \cong \BP^1$, $a_i \ge 2$ and $n/d=[a_1,a_2,\ldots,a_r]$ is the continued fraction expansion
\[
\frac{n}{d} = a_1 - {1\over\displaystyle a_2 - 
                                       {\strut 1\over\displaystyle \ddots - 
                                          {1\over\displaystyle a_{r-1} - 
                                            {1\over\displaystyle a_r}}}}
\]

For both cases where $G$ is non-cyclic or cyclic, let $\wt{\varphi}$ be a $G$-equivariant \'etale endomorphism of $\A^2$ 
such that $\wt{\varphi}^{-1}(O)=O$. Then we have the following result, which explains the singular fibers of 
$p : V \to C$ in terms of the $G$-action.

\begin{lem}\label{Lemma 1.1}
Let $\wh{\A}^2$ be the blowing-up of $\A^2$ at the point $O$ and let $E$ be the exceptional curve. Then the \'etale 
endomorphism $\wt{\varphi}$ induces an \'etale endomorphism $\wh{\varphi} : \wh{\A}^2 \to \wh{\A}^2$ such that 
the restriction of $\wh{\varphi}$ to $E$ is an automorphism. Furthermore, $\wh{\varphi}$ is $G$-equivariant, 
and $G$ acts on $E$ via the natural action of $\PGL(2,\C)$ on $\BP^1$. Let $\sigma$ be the induced action of $G$ 
on $E \cong \BP^1$. Then $\sigma$ has the following properties.
\begin{enumerate}
\item[(1)]
If $G$ is non-cyclic, the action $(G,\sigma)$ has three nontrivial isotropy subgroups up to conjugacy which are cyclic 
groups of order $m_1,\\ m_2,m_3$, where $\{m_1,m_2,m_3\}$ is one of the Platonic triplets. Consequently, there exist 
three cyclic singularities on $E/G \hookrightarrow \wh{\A}^2/G$. The minimal resolution of these singularities gives rise 
to the part $S_0+T_1+T_2+T_3$, which is one of the connected components of the boundary divisor of the normal completion $V$ 
of $X^\circ$, where $S_0$ is the proper transform of the curve $E/G$.
\item[(2)]
If $G$ is cyclic and $d > 1$, there are two $G$-fixed points $Q_0, Q_\infty$ on $E$ such that the $G$-action is locally 
described as 
\[
\begin{array}{ll}
(x,t) \mapsto (\zeta x, \zeta^{d-1}t) & {\text at}~~Q_0 \\
(u,y) \mapsto (\zeta^{n-d+1}u, \zeta^d y) & {\text at}~~Q_\infty\,
\end{array}
\]
where $\{x,y\}$ is the system of coordinates on $\A^2$, $t=y/x$ and $u=x/y$. In the minimal resolution, if 
$\delta_1=\gcd(n,d-1) >1$, we consider the reduced pair $(n/\delta_1, (d-1)/\delta_1)$ instead of $(n,d-1)$. Similarly, 
we consider the reduced pair for $(n-d+1,d)$. Namely, writing $\xi=\zeta^d$ and $\xi^e=\zeta^{n-d+1}~~(0 < e < n)$, 
we consider the pair $(n/\delta_2, e/\delta_2)$ with $\delta_2=\gcd(n,e)$. Let $S_0$ be the proper transform of 
$E/G$ and $T_0, T_\infty$ the exceptional loci of the minimal resolution of the image points of $Q_0, Q_\infty$, respectively. 
Then there exists a normal smooth completion $V$ and a $\BP^1$-fibration $p : V \to C$ such that 
$X^\circ=V-(S_0+T_0+T_\infty)-(S_1+R_0+R_\infty)$, where:
\begin{enumerate}
\item[(i)]
$S_0$ and $S_1$ are sections of $p$. Furthermore, $(S_0^2) \le -1$. 
\item[(ii)]
There exist $(-1)$-curves $F_0, F_\infty$ such that $T_0+F_0+R_0$ and $T_\infty+F_\infty+R_\infty$ are degenerate $\BP^1$-fibers.
Furthermore, the multiplicities $m_0$ and $m_\infty$ of $F_0$ and $F_\infty$ are respectively $n/\delta_1$ and $n/\delta_2$.
\item[(iii)]
The parts $R_0$ and $R_\infty$ are uniquely determined by the condition that $F_0$ and $F_\infty$ are unique $(-1)$-curves 
in the respective degenerate fibers.
\end{enumerate}
\item[(3)]
If $G$ is cyclic and $d=1$, then the induced $G$-action on the curve $E$ is trivial. Then $X^\circ=V-(S_0+S_1)$, where $V$ 
is a Hirzebruch surface $\F_n$, $S_0$ is a minimal section and $S_1$ is a section disjoint from $S_0$ 
\end{enumerate}
\end{lem}

\Proof
The endomorphism $\wt{\varphi}$ is defined by two polynomials $\varphi^\ast(x)=f(x,y)$ and $\varphi^\ast(y)=g(x,y)$, where 
\begin{eqnarray*}
&& f(x,y)=ax+by+\{\text{terms of degree higher than $1$}\} \\
&& g(x,y)=cx+dy+\{\text{terms of degree higher than $1$}\}\ .
\end{eqnarray*}
Let $t=y/x$. Then $\left.(g(x,y)/f(x,y))\right\vert_E=(c+dt)/(a+bt)$, where the Jacobian determinant 
$J((f,g)/(x,y))=ad-bc$ is a nonzero constant. Hence the endomorphism $\wt{\varphi}$ extends to an endomorphism 
$\wh{\varphi}$ of $\wh{\A}^2$. Since its restriction to $E$ is an automorphism, $\wh{\varphi}$ is an \'etale endomorphism. 
For this end, one can use purity of the branch locus to $\wh{\varphi}: \wh{\A}^2 \to \wh{\A}^2$. Furthermore, the ring 
endomorphism $\wt{\varphi}^\ast$ of $\C[x,y]$ is $G$-equivariant, its extension $\wh{\varphi}^\ast$ to the quotient field 
$\C(x,y)$ is also $G$-equivariant. Hence $\wh{\varphi}$ is $G$-equivariant.

Now we discuss separately the cases where $G$ is non-cyclic and $G$ is cyclic. If $G$ is non-cyclic, it is known that there 
exists a Galois covering $q : \wt{C} \to C$ with both $\wt{C}$ and $C$ isomorphic to $\BP^1$. The Galois group is $G$, and 
its action has the property stated in the assertion. See \cite[Lemma 4.3]{MM} for the references. In fact, the $G$-action $\sigma$ 
on $E$ is via $\PGL(2,\C)$, and it is known (see \cite{B}) that the covering $E \to E/G$ is exactly this Galois covering. Thence 
follows the assertions in (1).

Suppose that $G$ is cyclic and $d >1$. Most assertions are verified straightforwardly. If $\delta_1 > 1$, the subgroup $H=
\langle \zeta^{n/\delta_1}\rangle$ acts on the $(x,t)$-coordinates as pseudo-reflections which makes the variable $t$ invariant. 
Hence, by a theorem of Chevalley \cite{Chevalley}, the invariant subring is $\C[x,t]^H=\C[x',t]$ with $x'=x^{\delta_1}$ 
for which the induced $G$-action is given by ${}^{\zeta'}(x',t)=({\zeta'}x', {\zeta'}^{(d-1)/\delta_1}t)$, where 
$\zeta'=\zeta^{\delta_1}$ is a primitive $(n/\delta_1)$th root of $1$. Since the points $Q_0$ and $Q_\infty$ 
have cyclic singularities of the given type, the exceptional loci of these singularities connected by the component $S_0$ are 
obtaiend from the exceptional locus of the cyclic singularity of $\A^2/G$ by suitable subdivisional blowing-ups. Namely, these loci 
together with $S_0$ are contracted to a cyclic singularity of $\A^2/G$. Hence $(S_0^2) \le -1$. But it is not necessarily 
the case that $(S_0^2)=-1$. For example, if $n=6$ and $d=5$, $(S_0^2)=-2$. 

Suppose that $G$ is cyclic and $d=1$. Namely, the $G$-action on $\A^2$ is given by $(x,y) \mapsto (\zeta x,\zeta y)$. Then 
the action $\sigma$ on $E$ is trivial. Hence the quotient morphism $\wh{\pi} : \wh{\A}^2 \to \wh{\A^2}/G$ is totally ramifying 
on $E$. This implies that $\wh{\A}^2/G$ is smooth and $\wh{\pi}^\ast(S_0)=nE$. In fact, each orbit of the $G_m$-action 
induced by the standard $G_m$-action on $\A^2$, ${}^t(x,y)=(tx,ty)$, meets $E$ in one point transversally if taken the 
closure in $\wh{\A}^2$, where $E$ is totally ramifying for the morphism $\wh{\pi}$. Hence $\wh{\A}^2/G$ is an $\A^1$-bundle 
over the curve $S_0=\wh{\pi}(E)$. Hence $\wh{\A}^2/G$ is smooth. Since $n(S_0^2)=n^2(E^2)=-n^2$, it follows that 
$(S_0^2)=-n$ and $V$ is isomorphic to a Hirzebruch surface $\F_n$. 
\QED

Both in the non-cyclic and cyclic cases, we call $V$ the {\em standard} comletion of $X^\circ$, $p : V \to C$ the {\em standard} 
$\BP^1$-fibration and $p^\circ : X^\circ \to C$ the {\em standard} $\A^1_\ast$-fibration. Note that $p$ is determined by 
the $G_m$-action on $\A^2/G$ which is inherited from the natural $G_m$-action on $\A^2$ as the $G$-action on $\A^2$ is linear. 
Let $D=V-X^\circ$, which is a reduced effective divisor with simple normal crossings. The divisor $D$ has two connected components 
$D_0$ and $D_\infty$, where $D_0$ is the exceptional locus of resolution of the singular point $P_0$. We call $D+K_V$ the 
{\em logarithmic canonical divisor} (log canonical divisor, for short). 

\begin{lem}\label{Lemma 1.2}
The logarithmic canonical divisor $D+K_V$ is given by $\Sigma_i(\ell_i-F_i)-2\ell$, where $\ell$ is a general fiber of the $\BP^1$-fibration 
$p$ and the sum runs over all degenerate fibers of $p$. Hence we have
\begin{eqnarray*}
D+K_V &\sim& \left\{ \begin{array}{ll} \ell-(F_1+F_2+F_3) & \quad\text{\rm Case $G$ is non-cyclic} \\
                    -F_0-F_\infty & \quad \text{\rm Case $G$ is cyclic and $d > 1$} \\
                    -2\ell & \quad \text{\em Case $G$ is cyclic and $d =1$}
                    \end{array} \right.
\end{eqnarray*}               
\end{lem}
\Proof
The formula is obtained by the same computations as in \cite[p. 81]{M2}. 
\QED

\section{A quasi-\'etale endomorphism of $\A^2/G$}

Let $X=\A^2/G$ and let $\varphi : X_u \to X_l$ be a quasi-\'etale endomorphism. Let $\varphi^\circ : X^\circ_u \to 
X^\circ_l$ be the induced \'etale endomorphism. We first prove the following result which is a strengthening of 
\cite[Lemma 2.2]{MM}.

\begin{lem}\label{Lemma 2.1}\label{Lemma 2.1}
Let $X$ and $\varphi$ be the same as above. Then the following assertions hold.
\begin{enumerate}
\item[(1)]
There exists an \'etale endomorphism $\wt{\varphi} : \A^2 \to \A^2$ such that  $\wt{\varphi}^{-1}(O)=O$ and 
$\pi\cdot\wt{\varphi}=\varphi\cdot\pi$ for the quotient morphism $\pi : \A^2 \to X$. Furthermore, $\wt{\varphi}(gx)=
\chi(g)\wt{\varphi}(x)$ for a group endomorphism $\chi$ of $G$.
\item[(2)]
$\chi$ is an automorphism.
\item[(3)]
Replacing $\varphi$ by $\varphi^N$ for some $N > 0$, we may assume that $\chi=\id$.
\item[(4)]
Suppose $\varphi$ is birational. Then $\varphi$ is an automorphism.
\end{enumerate}
\end{lem}

\Proof
The assertion (1) is proved in the same fashion as in \cite[ibid.]{MM}. The argument therein works also in the case $G$ is cyclic. 
For the assertion (2), let $K=\Ker\chi$. Then $\wt{\varphi}(gx)=
\wt{\varphi}(x)$ for $g \in K$ and $x \in \A^2$. Then $\wt{\varphi}$ factors as 
\[
\wt{\varphi}^\circ : \A^2\setminus\{O\} \st{\pi_K^\circ} (\A^2/K)\setminus\{O\} \st{\tau^\circ} \A^2\setminus\{O\}\ .
\]
Since $\wt{\varphi}^\circ$ and $\pi_K^\circ$ are \'etale, so is $\tau^\circ$. On the other hand, the above three morphisms 
extend as
\[
\wt{\varphi} : \A^2 \st{\pi_K} \A^2/K \st{\tau} \A^2\ \ .
\]
Since $\wt{\varphi}$ is \'etale at the point $O$, $\wh{\SO}_{\A^2,\wt{\varphi}(O)} \cong \wh{\SO}_{\A^2,O}$. 
Let $Q_0=\pi_K(O)$, whence $\tau(Q_0)=O$. Then we have injections
\[
\SO_{\A^2,\wt{\varphi}(O)} \stackrel{\tau^\ast}{\hookrightarrow} \SO_{\A^2/K,Q_0} \hookrightarrow \SO_{\A^2,O},
\]
which induces the injections 
\[
\wh{\SO}_{\A^2,\wt{\varphi}(O)} \hookrightarrow \wh{\SO}_{\A^2/K,Q_0} \hookrightarrow \wh{\SO}_{\A^2,O}.
\]
In fact, $\wh{\SO}_{\A^2/K,Q_0} \to \wh{\SO}_{\A^2,O}$ is injective since $\wh{\SO}_{\A^2/K,Q_0}$ is the $K$-invariant 
subring of $\wh{\SO}_{\A^2,O}=\C[[x,y]]$. Since the composite 
\[
\wh{\SO}_{\A^2,\wt{\varphi}(O)} \to \wh{\SO}_{\A^2/K,Q_0} \to \wh{\SO}_{\A^2,O}
\]
is surjective, so is $\wh{\SO}_{\A^2/K,Q_0} \to \wh{\SO}_{\A^2,O}$. Hence $\wh{\SO}_{\A^2/K,Q_0} \cong \wh{\SO}_{\A^2,O}$. 
Hence $Q_0$ is a smooth point of $\A^2/K$. Since $K$ is small, this implies that $K=\{1\}$.
 
(3)~ We have $\wt{\varphi}^r(gx)=\chi^r(g)\wt{\varphi}^r(x)$. Since $G$ is a finite group, so is $\Aut(G)$. Hence $\chi^N=\id$ 
for some positive integer $N$. Then $\wt{\varphi}^N(gx)=g\wt{\varphi}^N(x)$. 

(4)~Replacing $\varphi$ by $\varphi^N$ for $N > 0$, by (3) above we may assume that $\varphi$ lifts to a $G$-equivariant 
\'etale endomorphism $\wt{\varphi}$ of $\A^2$. It is easy to show that $\wt{\varphi}$ is birational if so is 
$\varphi$. It is also well-known that $\wt{\varphi}$ is then an automorphism by Zariski's main theorem and a theorem of Ax-Grothendieck 
\cite{Ax,EGA}. Hence $\varphi : X \to X$ is an automorphism. 
\QED

We prove the following result. 

\begin{lem}\label{Lemma 2.2}
Let $X$ and $\varphi$ be the same as above. Suppose that $\varphi$ is induced by a $G$-equivariant \'etale endomorphism 
$\wt{\varphi} :\A^2 \to \A^2$ so that $\pi\circ\wt{\varphi}=\varphi\circ\pi$, where $\pi : \A^2 \to X$ is the quotient morphism. 
Then, after replacing $\varphi$ by its suitable power, $\varphi$ extends to an \'etale endomorphism $\wh{\varphi} : \wh{X} \to \wh{X}$ 
such that $\wh{\varphi}^{-1}(\wh{E})=\wh{E}$. Furthermore, $\wh{\varphi}\vert_{\wh{E}}$ is an automorphism.
\end{lem}

\Proof
Let $X$ be a complex normal algebraic surface with a quotient singular point $P_0$ and a quasi-\'etale endomorphism $\varphi$ 
of $X$ fixing the point $P_0$. We take $U$ to be a small analytic open neighborhood of $P_0$ such that $U$ is realized 
as a $G$-quotient of a small open disc $\wt{U}$ of $\A^2$ around the origin $O$, where $G=\pi_1(U\setminus\{P_0\})$. Let 
$B=\wh{\SO}_{U,P_0}$ and let $A=\wh{\SO}_{\wh{U},O}$. Let $\sigma : \Bl_{P_0}(X)\to X$ be the scheme-theoretic blowing-up of 
the point $P_0$. Namely, $\Bl_{P_0}(X)=\Proj(\oplus_{i \ge 0}\gm_{P_0}^i)$, where $\gm_{P_0}$ denotes, by abuse of notation, 
the sheaf of ideals defining the point $P_0 \in X$. The quasi-\'etale endomorphism $\varphi : X \to X$ induces the inclusion 
$\varphi^*(\gm_{P_0}) \subseteq \gm_{P_0}$. Let $M$ be the maximal ideal of $A$ and let $\gm=M\cap B$. We claim that $B=A^G$, 
$\gm$ is the maximal ideal of $B$ and $\varphi^*(\gm)=\gm$, where $\varphi^* : B \to B$ denotes, by abuse of notation, the ring 
homomorphism of complete local rings induced by $\varphi$. 

In fact, it is clear that $B=A^G$ and $\gm$ is the maximal ideal of $B$. After replacing $\varphi$ by its power if necessary, 
$\varphi|_U$ lifts up to a $G$-equivariant \'etale endomorphism $\wt{\varphi}$ of $\wt{U}$ such that $\wt{\varphi}^*(A)=A$ 
and $\wt{\varphi}^*(M)=M$, where we denote by the same letter the homomorphism of complete local rings induced by $\wt{\varphi}^*$ 
 (see Lemma \ref{Lemma 2.1}). It is clear that $\varphi^*(\gm) \subseteq \gm$ as $\wt{\varphi}^*(M)=M$ and $\wt{\varphi}^*(B)=B$. 
We show that $\varphi^*(\gm)=\gm$. For $b \in \gm$, there exists $a\in M$ such that $\wt{\varphi}^*(a)=b$ because $\wt{\varphi}^*(M)=M$. 
Set $b'=(\sum_{g\in G}{}^ga)/|G|$, which is an element of $\gm$ because $G$ stabilizes $M$. Then $\wt{\varphi}^*(b')=
(\sum_{g\in G}{}^g\wt{\varphi}^*(a)/|G|=(\sum_{g \in G}{}^gb)/|G|=b$. Hence $\wt{\varphi}^*(\gm)=\gm$. So, $\varphi^*(\gm)=\gm$.

Then the endomorphism $\varphi : X \to X$ lifts uniquely to an endomorphism $\psi : \Bl_{P_0}(X) \to \Bl_{P_0}(X)$ such 
that $\sigma\circ\psi=\varphi\circ\sigma$. Let $E$ be the exceptional locus of $\sigma$. Then $\psi(E)=E$. In fact, since 
$\varphi^*(\gm)=\gm$, each generic point of $E$ is in the image of $\psi|_E$. Since $E$ is proper, $\psi|_E$ is surjective. Replacing 
$\varphi$ by its suitable power if necessary, we may assume that $\psi|_E$ sends each irreucible component of $E$ to itself. Since 
$\psi|_E$ induces a separable algebraic extension of the function field of the irreducible exceptional component, $\psi$ is \'etale 
generically on each component of $E$ and hence \'etale on $E$ by purity of branch loci possibly except for the new singular points 
on $E$ of $\Bl_{P_0}(X)$.

Now revert to the notations in the statement. Since $\wt{\varphi} : \A^2 \to \A^2$ is \'etale at $O$ and $\varphi^{-1}(O)=O$, 
the above observation shows that $\varphi^*(\gm_{P_0})=\gm_{P_0}$ and $\varphi$ lifts uniquely to a quasi-\'etale  endomorphism 
$\psi : \Bl_{P_0}(X) \to \Bl_{P_0}(X)$.  Replacing $\varphi$ by its power, we may assume that $\psi$ fixes each of the singular 
points on $\Bl_{P_0}(X)$.

Let $\wh{\sigma} : \wh{X} \to X$ be the minimal resolution of singularity with exceptional locus $\wh{E}$. 
By Artin \cite[Proof of Theorem 4]{Artin}, $\gm_{P_0}\cdot\SO_{\wh{X}}=\SO_{\wh{X}}(-Z)$, where $Z$ is the fundamental cycle of 
$\wh{E}$. Hence $\wh{\sigma}$ factors through  the blowing-up $\sigma : \Bl_{P_0}(X) \to X$ as  
\[
\wh{\sigma} : \wh{X} \st{\tau_1} \Bl_{P_0}(X) \st{\sigma} X.
\]
If there exist new singular points on $\Bl_{P_0}(X)$ lying on the exceptional locus $E$ of $\sigma$, by the above factorization 
of $\wh{\sigma}$, each of the new singular point is obtained by contracting a partial subgraph of the exceptional locus $\wh{E}$ 
on $\wh{X}$, where we use the fact that $\Bl_{P_0}(X)$ is a normal surface (see the references below). Hence it is again 
a quotient singular point.
 
By Tyurina \cite[Theorem 1]{Tyurina} or by Lipman \cite[Proposition 8.1 and Theorem 4.1]{Lipman}, $\Bl_{P_0}(X)$ is a normal 
algebraic surface and $\wh{X}$ is obtained by a finite succession of blowing-ups of the above type. Hence $\varphi$ extends 
to an \'etale endomorphism $\wh{\varphi} : \wh{X} \to \wh{X}$ such that $\wh{\varphi}^{-1}(\wh{E})=\wh{E}$. 

As for the last assertion, note that the restriction $\wh{\varphi}\vert_{\wh{E}}$ is a proper \'etale endomorphism since 
$\wh{E}$ is a proper scheme. Hence it is a topological covering of the simply-connected space $\wh{E}$. This implies that 
$\wh{\varphi}\vert_{\wh{E}}$ is an automorphism. 
\QED

\begin{remark}\label{Remark 2.3}{\em 
In order to prove Theorem stated in the introduction, we use Lemma \ref{Lemma 2.2} only in the case where $X$ has a cyclic 
quotient singularity of type $(n,1)$ (or type $\frac{1}{n}(1,1)$ in the current notation). In this case, Lemma 2.2 is proved 
in a more transparent way. 

Let $\sigma_1 : \wh{\A}^2 \to \A^2$ be the blowing-up of the origin $O$ and let $E_1$ be the exceptional curve. By Lemma 
\ref{Lemma 1.1}, $\wt{\varphi}$ induces a $G$-equivariant \'etale endomorphism $\wh{\varphi} : \wh{\A}^2 \to \wh{\A}^2$ 
such that $\wh{\varphi}^{-1}(E_1)=E_1$, $E_1$ is $G$-stable and $\wh{\varphi}\vert_{E_1}$ is an automorphism. This assertion holds 
because $O$ is a $G$-fixed point and $\wt{\varphi}^*$ induces a local ring automorphism of $\wh{\SO}_{\A^2,O}$. Let $X_1=\wh{\A}^2/G$ 
and let $\ol{E}_1=E_1/G$. Since $\wh{\varphi}$ is $G$-equivariant and $\wh{\varphi}^{-1}(E_1)=E_1$, $\wh{\varphi}$ induces 
an endomorphism $\varphi_1 : X_1 \to X_1$ such that $q_1\circ \wh{\varphi}=\varphi_1\circ q_1$ and $\varphi_1^{-1}(\ol{E}_1)=\ol{E}_1$, 
where $q_1 : \wh{\A}^2 \to X_1$ is the quotient morphism. 

Since the $G$-action has type $(n,1)$, the induced $G$-action on $E_1$ is trivial. Then $X_1$ is smooth, $q_1$ is totally 
ramified over $\ol{E}_1$ and $\varphi_1$ is \'etale on $\ol{E}_1$. Let $X$ and $\varphi$ be the same as in Lemma \ref{Lemma 2.2}. 
Then the endomorphism $\wh{\varphi}$ induces an automorphism on the exceptional locus $\wh{E}$. This implies that the given 
quasi-\'etale endomorphism $\varphi$ is \'etale at the singular point $P_0$ as well. \QED}
\end{remark}

We prove the following result, which will be used in the subsequent arguments.

\begin{thm}\label{Theorem 2.4}
Let $\varphi : X \to X$ be a quasi-\'etale endomorphism of $X=\A^2/G$ for a small finite subgroup $G$ of $\GL(2,\C)$. Suppose that 
$\varphi$ preserves the standard $\A^1_\ast$-fibration $p^\circ$. Namely there exists an endomorphism $\beta : C \to C$ such that 
$p^\circ\cdot(\varphi\vert_{X^\circ})=\beta\cdot p^\circ$. Then $\varphi$ is an automorphism.
\end{thm}
\Proof
The following simple proof was suggested by one of the referees. 

By Lemma \ref{Lemma 2.2}, we assume after replacing $\varphi$ by its suitable power that $\varphi$ extends to an \'etale 
endomorphism $\wh{\varphi}$ of $\wh{X}$ which induces an automorphism on the curve $S_0$ with the notation in section one. Since 
$S_0$ is identified with $C$ via $p|_{S_0}$, $\wh{\varphi}|_{S_0}$ is identified with $\beta$. Hence $\beta$ is an automorphism. 

For a general fiber $\ell$ of the $\BP^1$-fibration $p$, $\ell\setminus(\ell\cap S_1)$ is isomorphic to $\A^1$, and $\wh{\varphi}$ 
restricted to $\ell\setminus(\ell\cap S_1)$ induces an \'etale endomorphism from $\A^1$ to $\A^1$, which must be an isomorphism by 
the Jacobian conjecture in dimension one. Hence $\wh{\varphi}$ and $\varphi$ are birational. Then $\varphi$ is an automorphism 
by Lemma \ref{Lemma 2.1}. 
\QED

\section{When does $\varphi$ preserve the standard $\A^1_\ast$-fibration?}

Let $\nu : Y \to X_l$ be the normalization morphism with $Y$ the 
normalization of $X_l$ in the function field of $X_u$. Then, by Zariski's main theorem, there exists an open immersion 
$\iota : X_u \to Y$ such that $\nu\cdot \iota=\varphi$. We identify $X_u$ with the image $\iota(X_u)$. By \cite[Th. 1.1]{GM}, 
we know the following result.

\begin{lem}\label{Lemma 3.1}
$Y-X_u$ consists of a disjoint union of irreducible curves $\coprod_{i=1}^r A_i$ such that $A_i$ is isomorphic to $\A^1$ 
and, for each $i$, there exists at most one cyclic quotient singularity $P_i$ lying on $A_i$. The exceptional locus $B_i$ of the 
cyclic singularity at $P_i$ meets the proper transform of $A_i$ at the end component of $B_i$.
\end{lem}

\begin{remark}\label{Remark 3.2}{\em
We consider the set $(Y-X_u)\cap \nu^{-1}(P_0)$. 

(1)~Assume first that $(Y-X_u)\cap \nu^{-1}(P_0) = \emptyset$. Then $\varphi : X_u \to X_l$ is proper over $P_0$. This implies that 
$\varphi^* : \SO_{X_l,P_0} \to \SO_{X_u,P_0}$ is an isomorphism. In fact, since $\nu^{-1}(P_0)=P_0$ by the assumption, 
$\SO_{X_u,P_0}$ is a finite $\SO_{X_l,P_0}$-module. Since $\varphi$ is induced by a $G$-equivariant endomorphism 
$\wt{\varphi}$ of $\A^2$, by Lemma \ref{Lemma 2.2}, we have 
\[
\SO_{X_u,P_0}\otimes_{\SO_{X_l,P_0}}\wh{\SO}_{X_l,P_0}\cong \wh{\SO}_{X_u,P_0} \cong \wh{\SO}_{X_l,P_0} 
\cong \SO_{X_l,P_0}\otimes_{\SO_{X_l,P_0}}\wh{\SO}_{X_l,P_0}.
\]
Since $\wh{\SO}_{X_l,P_0}$ is faithfully flat over $\SO_{X_l,P_0}$, it follows that $\SO_{X_u,P_0}\cong \SO_{X_l,P_0}$. 
This implies that $\varphi$ is birational. Then $\varphi$ is an automorphism by Lemma \ref{Lemma 2.1}, (4).

(2)~Suppose that the image by $\nu$ of an irreducible component of $Y-X_u$ passes through the point $P_0$. Namely, we suppose that 
$\nu(P)=P_0$ for a point $P$ of the component $A_i$, say $i=1$. Let $(U/H,P)$ be a local analytic neighborhood of $P$ in $Y$ 
such that $(U/H,P)$ is mapped to a local analytic neighborhood $(U_0/G,P_0)$ of $X$ at the point $P_0$ via a finite 
holomorphic map $\nu_1$ of degree $n_1$ induced by $\nu$, where $U$ and $U_0$ are local analytic neighborhoods of $\A^2$ at 
the origin $O$ and $H=\pi_1(U/H-\{P\})$ is a cyclic group which is trivial if $P$ is a smooth point. Denoting by the superscript 
$\circ$ the germs punctured at the origin or the singular point, we have the following commutative diagram.

\raisebox{-60mm}{
\begin{picture}(145,70)(10,0)
\setlength{\unitlength}{0.9mm}
\put(40,60){$U^\circ$}
\put(35,55){\vector(-3,-2){15}}
\put(45,55){\vector(3,-2){15}}
\put(0,40){$U_0^\circ\times_{(U_0/G)^\circ}(U/H)^\circ$}
\put(40,40){\vector(1,0){15}}
\put(60,40){$(U/H)^\circ$}
\put(22,33){\vector(0,-1){10}}
\put(65,33){\vector(0,-1){10}}
\put(20,15){$U_0^\circ$}
\put(33,16){\vector(1,0){20}}
\put(60,16){$(U_0/G)^\circ$}
\put(55,55){$H$}
\put(45,43){$G$}
\put(42,8){$G$}
\put(60,28){$\nu_1$}
\put(15,28){${\rm pr}$}
\put(80,16){$,$}

\put(85,45){$(U/H_0)^\circ$}
\put(103,45){\vector(1,0){10}}
\put(115,45){$(U/H)^\circ$}
\put(90,37){\vector(0,-1){10}}
\put(95,30){$\mu$}
\put(120,37){\vector(0,-1){10}}
\put(123,30){$\nu_1$}
\put(87,18){$U_0^\circ$}
\put(100,20){\vector(1,0){13}}
\put(117,18){$(U_0/G)^\circ$}
\end{picture}}

In the left picture, the lower square is a fiber product. Hence the upper horizontal arrow indicates an \'etale finite  
covering. Hence $U^\circ_0\times_{(U_0/G)^\circ}(U/H)^\circ$ is a direct sum of connected components $(U^\circ/H_i)^\circ$ 
for a subgroup $H_i$ of $H$. Since $G$ acts transitively on the connected components, $H_i$ is equal to one and the same subgroup 
$H_0$ of $H$. The left vertical arrow of the square is a finite map. Hence, by extracting these data, we obtain the right 
picture, where $\mu$ denotes the restriction of the projection map to one connected component of 
$U^\circ_0\times_{(U_0/G)^\circ}(U/H)^\circ$. Hence we have $[H:H_0]n_1=m|G|$, where $n_1=\deg \nu_1$ and $m=\deg \mu$. One 
connected component $(U/H_0)^\circ$ has degree $[H:H_0]$ over $(U/H)^\circ$, and its translate $g(U/H_0)^\circ \cong (U/H_0)^\circ$ has 
the same degree $[H:H_0]$ over $(U/H)^\circ$. Hence $[H:H_0]$ divides $|G|$ by the count of degree of the upper horizontal morphism of 
the left picture, and consequently, $m$ divides $n_1$. In particular, if $[H:H_0]=1$, e.g., if $P$ is smooth, then $n_1$ is 
a multiple of $|G|$. 

(3)~If the morphism $\nu_1$ is unramified on the component $A_1$, then $\mu$ is an isomorphism in the above diagram because 
the restriction to $A_1$ of the projection morphism is finite and unramified. This implies that $m=1$ and the universal covering map 
$U_0^\circ \to (U_0/G)^\circ$ is factored as
\[
U_0^\circ \cong (U/H_0)^\circ \to (U/H)^\circ \st{\nu} (U_0/G)^\circ\ .
\]
Hence $U$ is homotopic to $U_0$ and  $H$ is a subgroup of $G$. 
} \QED
\end{remark} 

Hereafter we assume without loss of generality that $\varphi$ lifts to a $G$-equivariant endomorphism $\wt{\varphi}$ of $\A^2$. 
Let $\wh{Y}$ be a minimal resolution of singularities on $Y$. Then $\wh{Y}$ contains $\wh{X}_u$ as an open set. By blowing up 
the points on $\wh{Y}-\wh{X}_u$, we find a minimal sequence of blowing-ups $\alpha : \wt{Y} \to \wh{Y}$ such that the canonical 
rational mapping 
\[
\wh{Y} \dto \wh{X}_u \st{\wh{\varphi}} \wh{X}_l
\] 
composed with $\alpha$ is a morphism $\wt{\nu} : \wt{Y} \to \wh{X}_l$. The above rational mapping is determined by the 
inclusion of the function fields
\[
\C(\wh{X}_l) \stackrel{\wh{\varphi}^*}{\hookrightarrow} \C(\wh{X}_u) =\C(\wh{Y}),
\]
which might not be a morphism if there exist points on $Y\setminus X_u$ mapped by $\nu$ to the point $P_0$. Such points 
on $Y$ must be blown up. Then $\wt{\nu}\vert_{\wh{X}_u}=\wh{\varphi}$  by Lemma \ref{Lemma 2.2}. Let $W$ be 
a log smooth completion of $\wt{Y}$, i.e., a smooth projective surface such that $W-\wt{Y}$ is a divisor with simple
 normal crossings, such that the birational mapping $W \dto V_u$ is a morphism which is a composite of blowing-ups and the morphism 
$\wt{\nu}$ extends to a projective morphism $\Phi : W \to V_l$ which makes the following diagram commutative
\[
\CD
\wt{Y} @>{}>> W \\
@V\wt{\nu} VV   @VV\Phi V  \\
\wh{X}_l @>{}>> V_l 
\endCD
\]
Here $\wh{X}_l$ is the minimal resolution of singularity of $X_l$, and $V_l$ is the standard completion of 
$X^\circ_l$. Since $\wt{Y}$ contains $\wh{X}_u$, $W$ is obtained from the standard completion $V_u$ of $X^\circ_u$ by 
blowing up points on the boundary divisor $D_u (= D)$. Let $\rho : W \to V_u$ be the birational morphism which is a 
composite of blowing-ups. The irreducible component $A_i$ of $Y-X_u$ together with the exceptional locus $B_i$ of the 
possible cyclic singularity on $A_i$ is obtained from a point on an irreducible component of $D$ as exhibited in 

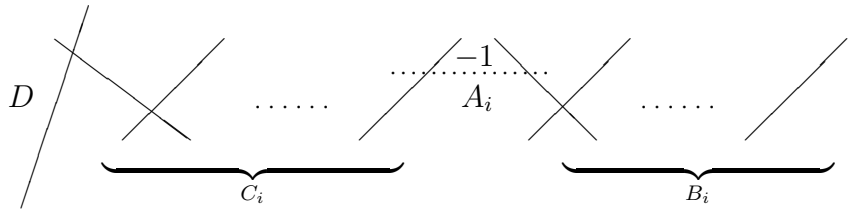
\raisebox{-35mm}{
\begin{picture}(130,40)(0,0)
\unitlength=0.9mm

\put(5,5){\line(1,3){10}}
\put(10,30){\line(4,-3){20}}
\put(20,15){\line(1,1){15}}
\multiput(40,20)(2,0){6}{\circle*{0.3}}
\put(55,15){\line(1,1){15}}
\multiput(60,25)(1.5,0){16}{\circle*{0.2}}
\put(75,30){\line(1,-1){15}}
\put(80,15){\line(1,1){15}}
\multiput(97,20)(2,0){6}{\circle*{0.3}}
\put(112,15){\line(1,1){15}}
\put(3,20){$D$}
\put(69,26){$-1$}
\put(70,20){$A_i$}
\put(85,12){$\underbrace{\hspace{36mm}}_{B_i}$}
\put(17,12){$\underbrace{\hspace{40mm}}_{C_i}$}
\end{picture}}

\noindent
where $A_i$ is a $(-1)$-curve \footnote{$A_i$ is the last exceptional $(-1)$-curve in the blowing-up process and the part $B_i$ 
is the locus of exceptional curves of a minimal resolution of a cyclic quotient singularity. Hence $A_i$ stays as a $(-1)$-curve.}, 
$B_i$ is a linear chain not necessarily consisting of only $(-2)$-curves and the part $C_i$
connecting $D$ and $A_i$ is not necessarily a linear chain. The divisor $A_i+B_i+C_i$ is contracted to a smooth point on $D$ or 
on the divisor obtained from $D$ by subdivisional blowing-ups. In the midst of blowing-ups $\rho : W \to V_u$, the following 
five cases are possible. We denote by $\sigma$ an intermediate single blowing-up. The change of log canonical divisors under 
$\sigma$ is given as follows. 
\begin{enumerate}
\item[(1)]
If $\sigma$ is a blowing-up of a non-boundary point with the exceptional curve $E$ included in the new boundary divisor, 
then $D''+K_{V''}=\sigma^*(D'+K_{V'})+2E$, where $\sigma : V'' \to V'$ is the blowing-up, $D'$ is the boundary divisor of $V'$ and 
$D''=\sigma'(D')+E$ is the boundary divisor of $V''$. This case occurs if a point on $A_i$ is mapped to $P_0$ by $\nu$.
\item[(2)]
If $\sigma$ is a subdivisional blowing-up of the boundary divisor with the resulting exceptional curve $E$ included in the new 
boundary divisor, then $D''+K_{V''}=\sigma^\ast(D'+K_{V'})$, where $D''=\sigma'(D')+E$.
\item[(3)]
If $\sigma$ is a subdivisional blowing-up with $E$ not included in the boundary divisor, then $D''+K_{V''}=\sigma^\ast(D'+K_{V'})-E$.
\item[(4)]
If $\sigma$ is a sprouting blowing-up with $E$ not included into the boundary divisor, then $D''+K_{V''}= 
\sigma^\ast(D'+K_{V'})$.
\item[(5)]
If $\sigma$ is a sprouting blowing-up with $E$ included in the boundary divisor, then $D''+K_{V''}=\sigma^\ast(D'+K_{V'})+E$.
\end{enumerate}

Since $\wh{Y}-\wh{X}_u=\coprod_{i=1}^r(A_i+B_i)$, the cases (3) and (4) occur only for $E=A_i$. By abuse of notation, 
we here denote by the same $A_i$ the closure of $A_i$ in $\wh{Y}$ (later also in $\wh{V}, \wt{V}$ and $W$). More precisely, 
(3) occurs if $B_i \ne 0$ and (4) occurs if $B_i=0$, i.e., there is no singular point of $Y$ on $A_i$. Further, the rational morphism 
$\alpha : \wt{Y} \to \wh{Y}$ is a composite of blowing-ups with centers on $\wh{Y}-X_u$. 

We examine more precisely the change of log canonical divisors by $\rho : (W,\Delta) \to (V_u,D)$, where $\Delta=W-Y^\circ$ with 
the smooth locus $Y^\circ$ of $Y$. Decompose the morphism $\rho : W \to V_u$ as in the top row of the diagram below, where $\wh{V}$ 
(resp. $\wt{V}$) is a log smooth completion of $\wh{Y}$ (resp. $\wt{Y}$) such that $\wt{V}-\wt{Y} \cong \wh{V}-\wh{Y}$. 
Hence the blowing-ups appearing in the morphism $\rho_2$ are those appearing in the morphism $\alpha$, which are performed to 
eliminate the indeterminacy of the rational mapping $\wh{Y} \dto \wh{X}_u \st{\wh{\varphi}} \wh{X}_l$. Those centers are points 
on the $A_i$ mapped to the point $P_0$ 
of $X_u$ and their infinitely near points. Note that the $\BP^1$-fibration $p : V_u \to C$ extends to $\BP^1$-fibrations 
$\wh{p} : \wh{V} \to C$ and $\wt{p} : \wt{V} \to C$ such that the divisors $A_i+B_i$ are contained in the fibers of $\wh{p}$. 
Thus there exists a $\BP^1$-fibration $p_W : W \to C$ which is induced by the $\BP^1$-fibration $p : V_u \to C$, and each irreducible 
component of $\wh{Y}-X_u$ is contained in the fibers of $p_W$. This shows that $\wh{Y}-X_u$ (and $W-X_u$ as well) is a divisor 
with simple normal crossings. 

\raisebox{-60mm}{
\begin{picture}(100,70)(0,0)
\unitlength=0.9mm

\put(4,60){$\rho~~:$}
\put(15,60){$W$}
\put(20,60){\vector(1,0){15}}
\put(26,63){$\rho_3$}
\put(40,60){$\wt{V}$}
\put(45,60){\vector(1,0){15}}
\put(51,63){$\rho_2$}
\put(65,60){$\wh{V}$}
\put(70,60){\vector(1,0){15}}
\put(76,63){$\rho_1$}
\put(90,60){$V_u$}
\put(20,55){\vector(1,-1){40}}
\put(40,40){$\wt{Y}$}
\put(45,41){\vector(1,0){15}}
\put(51,44){$\alpha$}
\put(65,40){$\wh{Y}$}
\put(72,40){$\cdots\dto$}
\put(90,40){$\wh{X}_u$}
\put(40,50){$\bigcup$}
\put(65,50){$\bigcup$}
\put(90,50){$\bigcup$}
\put(93,35){\vector(0,-1){17}}
\put(47,35){\vector(2,-1){40}}
\put(90,10){$\wh{X}_l$}
\put(70,27){$\wt{\nu}$}
\put(65,10){$V_l$}
\put(80,10){$\supset$}
\put(35,30){$\Phi$}
\put(95,27){$\wh{\varphi}$}
\end{picture}}

We observe the changes of log canonical divisors according to the decomposition $\rho=\rho_1\circ\rho_2\circ\rho_3$.
\begin{enumerate}
\item[(I)]
Let $\wh{\Delta}=\wh{V}-Y^\circ$, which is considered as a reduced effective divisor. Then we have the following formula 
for the morphism $\rho_1 : \wh{V} \to V_u :$ 
\[
\wh{\Delta}+K_{\wh{V}}=\rho_1^*(D+K_{V_u})+\sum_j E_j-\sum_{i=1}^s A_i,
\]
where the $E_j$ exhaust the total transforms on $\wh{V}$ of the exceptional curves which arise from sprouting blowing-ups possibly after 
subdivisional blowing-ups of the divisor $D$ and are included in the boundary divisor and where the $A_i~(1 \le i \le s)$ exhaust the 
irreducible components of $Y-X_u$ carrying singular points of $Y$. Here $\sum_jE_j-\sum_{i=1}^sA_i$ is an effective divisor. 

In fact, the completion $\wh{V}$ of $\wh{Y}$ is obtained from $V_u$ by means of sprouting or subdivisional blowing-ups of the 
connected component $D_\infty$ of the boundary $D$ of $V_u$ and subsequent new boundaries, whereas the connected component $D_0$ of $D$ 
arising from the resolution of singularity $P_0$ of $X_u$ is untouched. Effect of the new exceptional curve in the log canonical divisor 
in each step of blowing-ups is considered in the above cases (2) $\sim$ (5). In particular, the component $A_i$ with $B_i \ne \emptyset$ 
is obtained from a sprouting blowing-up in an intermediate step followed by subdivisional blowing-ups. Then the total transform on $\wh{V}$ 
of the exceptional curve arising from the sprouting one contains the last exceptional $(-1)$ curve $A_i$. Hence $\sum E_j-\sum_{i=1}^sA_i$ is 
effective.

\item[(II)]
By Remark \ref{Remark 3.2}, (1), we may assume that there exists a point on $\wh{Y}-X_u$ which is mapped to the point $P_0$ on $X_l$. 
Since the singularity at $P_0$ is resolved on $\wh{V}_l$, we need the corresponding blowing-ups at the points of $\wh{Y}$ mapped to 
$P_0$. Hence the morphism $\rho_2 : \wt{V} \to \wh{V}$ coincides with the morphism $\alpha : \wt{Y} \to \wh{Y}$ on the open set $\wt{Y}$ 
and the identity on the boundary $\wt{V}-\wt{Y}$. The morphism $\alpha$ is a composite of blowing-ups with centers at points on $A_i+B_i$ 
(including the case $B_i=\emptyset$) and their infinitely near points. All resulting exceptional curves are contained in the boundary 
divisor $\Delta=W-Y^\circ$ when taken the proper transforms on $W$. The change of log canonical divisors is given by 
\[
\wt{\Delta}+K_{\wt{V}}=(\rho_1\rho_2)^*(D+K_{V_u})+\wt{\Gamma}_1+\wt{\Gamma}_2-\sum_{i=1}^sA'_i,
\]
where $\wt{\Gamma}_1=\rho_2^*(\sum_jE_j)$, $\wt{\Gamma}_2$ is a sum of the total transforms of the divisors $\delta_EE$ with $E$ ranging 
over all exceptional curves arising from $\alpha$ and $A'_i=\rho'_2(A_i)$ is the proper transform of $A_i$. The coefficient $\delta_E$ is 
equal to $2$ if the blowing-up is centered on $A_i$ with $B_i=\emptyset$ (by the case (1) above), it is equal to $1$ if the blowing-up 
is centered on $A_i\setminus(A_i\cap B_i)~(1 \le i \le s)$ (since $A_i$ produces $E$ once) or if the blowing-up is sprouting one centered 
on a boundary component which is not the proper transform of $A_i~(1 \le i \le s)$ (by the cases (2) and (5)), and it is $0$ otherwise. 
Here $\wt{\Gamma}_1+\wt{\Gamma}_2-\sum_{i=1}^sA'_i$ is an effective divisor.

\item[(III)]
The birational morphism $\rho_3 : W \to \wt{V}$ consists of blowing-ups with centers (including infinitely near centers) over 
$\wt{V}-\wt{Y}$. The log canonical divisor on $W$ is given by 
\[
\Delta+K_W=\rho^*(D+K_{V_u})+\Gamma_1+\Gamma_2+\Gamma_3-\sum_{i=1}^sA'_i,
\]
where $\rho=\rho_1\rho_2\rho_3$, $\Gamma_1=\rho_3^*(\wt{\Gamma}_1)$, $\Gamma_2=\rho_3^*(\wt{\Gamma}_2)$, $\Gamma_3$ is the sum of total 
transforms of sprouting blowing-ups centered over $\wt{V}-\wt{Y}$ and $A'_i$ is the proper transform of $A_i$. The divisor added 
to $\rho^*(D+K_{V_u})$ on the right hand side is effective. 
\end{enumerate}

Summarizing the above arguments, we obtain the following result.

\begin{lem}\label{Lemma 3.3}
Let $\Delta=W-Y^\circ$ with the smooth locus $Y^\circ$ of $Y$. Then we have
\[
\Delta+K_W=\rho^\ast(D+K_{V_u})+\Gamma-\sum_{i=1}^s A'_i\ ,
\]
where $\Gamma$ is an effective divisor supported by the exceptional divisor of the birational morphism $\rho : W \to V_u$ and 
$A'_i$ is the proper transform of $A_i$ with $B_i \ne \emptyset$. Further, $\Gamma \ge \sum_{i=1}^sA'_i$.
\end{lem}

The logarithmic ramification divisor formula due to Iitaka \cite{Iitaka} (see \cite[Lemma 1.11.1]{M3}) applied to 
the morphism $\Phi : (W,\Delta) \to (V_l,D)$ implies that
\[
\Delta+K_W=\Phi^\ast(D+K_{V_l})+R,
\]
where $R$ is an effective divisor supported by the union of curves contracted by $\Phi$ and the curves on which $\Phi$ is 
finite to the image but ramifying. We describe the ramification locus $R$. 

\begin{lem}\label{Lemma 3.4}
Let $V, W$ be smooth projective surfaces and let $D, \Delta$ be reduced effective divisors with simple normal crossings on 
$V, W$, respectively. Let $\Phi : W \to V$ be a dominant morphism such that $\Phi^{-1}(D)\subseteq \Delta$ and every divisor 
component in $\Delta\setminus\Phi^{-1}(D)$ is contracted by $\Phi$ \footnote{Hence it follows that $\Phi_*(\Delta)=D$}. 
Let $R$ be the logarithmic ramification divisor defined by the above formula. Then $R$ is supported by the union of curves $C$ 
on $W$ such that $C$ is contracted by $\Phi$ or $\Phi\vert_C : C \to \Phi(C)$ is ramifying. The curve $C$ has coefficient zero 
in $R$ if $C$ is a component of $\Delta$ which is not contracted by $\Phi$. If a component $C$ of $\Delta$ has coefficient zero 
in $R$ and $C$ is contracted by $\Phi$, then it is contracted to an intersection point of two irreducible components of $D$.
\end{lem}

\Proof
Suppose that $C$ is not contracted by $\Phi$. Let $\ol{C}=\Phi(C)$. Choose a general point $P$ on $C$ and let $Q=\Phi(P)$. 
Suppose that $C$ is not a component of $\Delta$. We choose a system of parameters $\{t,u\}$ of $V$ at $Q$ so that $\ol{C}$ 
is locally defined by $t=0$ and $u$ is a parameter along $\ol{C}$. Similarly, we choose a system of parameters $\{\xi,\eta\}$ 
of $W$ at $P$ so that $C$ is defined by $\xi=0$ and $\eta$ is a parameter along $C$. Then $\Phi^\ast(t) =f(\xi,\eta)$ and 
$\Phi^\ast(u)=g(\xi,\eta)$, where we view $f, g$ as complex functions holomorphic at $P$. We fix a nonzero, rational, 
differential $2$-form $\omega$ of $V$ and express it as $\omega=adt\wedge du$, where $a \in \C(V)$. Similarly, write 
$\Phi^\ast(\omega)=bd\xi\wedge d\eta$ with $b \in \C(W)$. The divisor $D+K_V$ (resp. $\Delta+K_W$) is determined near the point 
$Q$ (resp. $P$) by evaluating $b$ (resp. $a$) in terms of the normalized discrete valuation $v_C$ (resp. $v_{\ol{C}}$) 
associated to $C$ (resp. $\ol{C}$). By the assumption, $\ol{C}$ is not a component of $D$ as well. By a simple calculation, 
we obtain
\begin{eqnarray*}
b=Ja \qquad \text{with}\qquad J=J\left(\frac{f,g}{\xi,\eta}\right).
\end{eqnarray*}
Hence, near the point $P$, $\Delta+K_W$ differs from $\Phi^\ast(D+K_V)$ by $v_C(J)C$. Note that $\Phi\vert_C : C \to \ol{C}$ 
is unramified at $P$ if and only if $J(P) \ne 0$. In fact, if $v_C(f)=r$, then $v_C(J)=r-1$. 

If $C$ is a component of $\Delta$, then $\ol{C}$ is a component of $D$. Express $\omega$ and $\Phi^\ast(\omega)$ as 
\begin{eqnarray*}
\omega=a\frac{dt}{t}\wedge du \qquad \text{and} \qquad \Phi^\ast(\omega)=b\frac{d\xi}{\xi}\wedge d\eta. 
\end{eqnarray*}
Then we have 
\begin{eqnarray*}
b=\frac{a\xi}{f}J\left(\frac{f,g}{\xi,\eta}\right).
\end{eqnarray*}
Hence, near the point $P$, $\Delta+K_W$ differs from $\Phi^\ast(D+K_V)$ by $v_C(\xi J/f)C$ which is zero. In fact, if we write 
$f=\xi^mf_1$ with $m >0$ and $\xi \nmid f_1$ in $\SO_{W,P}$ then we have 
\begin{eqnarray*}
\frac{\xi J}{f}&=& \frac{m\xi^mf_1g_\eta+\xi^{m+1}J_1}{\xi^mf_1} =mg_\eta+\frac{\xi J_1}{f_1},
\end{eqnarray*}
where $J_1$ is the Jacobian of $f_1, g$ with respect to $\xi, \eta$ and $v_C(g_\eta)=0$ for a general point $P$ of $C$.
\svskip

Now suppose that $C$ is contracted by $\Phi$. Then, choosing a system of parameters $\{t, u\}$ of $V$ at the point $Q=\Phi(C)$, 
one can write $\Phi^\ast(t)=\xi^r f_1(\xi,\eta)$ and $\Phi^\ast(u)=\xi^s g_1(\xi,\eta)$ with $\xi \nmid f_1g_1,~r>0$ and $s>0$. Hence 
we have   
\begin{eqnarray*}
J\left(\frac{f,g}{\xi,\eta}\right)&=&\xi^{r+s-1}\left(rf_1\frac{\partial g_1}{\partial \eta}-sg_1\frac{\partial f_1}{\partial \eta}\right)
+\xi^{r+s}J\left(\frac{f_1,g_1}{\xi,\eta}\right) \\
&=& \xi^{r+s-1}I \quad \text{with} \quad I=rf_1\frac{\partial g_1}{\partial \eta}-sg_1\frac{\partial f_1}{\partial \eta}+
\xi J\left(\frac{f_1,g_1}{\xi,\eta}\right).
\end{eqnarray*}
This implies that
\begin{eqnarray*} 
\omega=adt\wedge du =a\xi^{r+s-1}I d\xi\wedge d\eta.
\end{eqnarray*} 
If $C$ is not a component of $\Delta$, then $\Delta+K_W$ differs from $\Phi^\ast(D+K_V)$ by $v_C(\xi^{r+s-1}I)C$ which has a positive 
coefficient. If $C$ is a component of $\Delta$, three cases are possible. Namely, the point $Q=\Phi(C)$ lies in $V\setminus D$, lies on 
a single irreducible component $D_1$ or is the intersection point $D_1\cap D_2$. In the first case, 
\[
\omega=adt\wedge du =a\xi^{r+s}I\frac{d\xi}{\xi}\wedge d\eta
\]
and the coefficient of $C$ in $(\Delta+K_W)-\Phi^*(D+K_V)$ is $v_C(\xi^{r+s}I)$ which is positive. In the second case, we choose 
a system of local parameters $\{t,u\}$ so that $D_1$ is defined by $t=0$ and $u$ is a parameter along $D_1$. In the third case, 
we choose $\{t,u\}$ so that $t$ and $u$ are respectively parameters along the curves $D_2$ and $D_1$. In the second case, we have
\begin{eqnarray*}
\omega=a\frac{dt}{t}\wedge du=\frac{a\xi^s I}{f_1}\frac{d\xi}{\xi}\wedge d\eta,
\end{eqnarray*}
where $f_1 \ne 0$ on a general point of $C$. In the third case, we have 
\begin{eqnarray*}
\omega=a\frac{dt}{t}\wedge \frac{du}{u}=\frac{aI}{f_1g_1}\frac{d\xi}{\xi}\wedge d\eta.
\end{eqnarray*}
In the second case, the component $C$ appears with a positive coefficient in the difference $(\Delta+K_W)-\Phi^\ast(D+K_V)$. 
In the third case, the coefficient is non-negative with $C$ possibly not appearing in the difference.
\QED

Comparing the logarithmic ramification divisor formula with the one in Lemma \ref{Lemma 3.3},
we have the following formula.
\begin{eqnarray}
\Phi^\ast(D+K_{V_l})=\rho^\ast(D+K_{V_u})+\Gamma-\sum_{i=1}^sA'_i-R. 
\end{eqnarray}
Let $H$ be the irreducible boundary component of $W-\wh{X}_u$ which is the proper transform of the section $S_1$ of the standard 
$\BP^1$-fibration $p : V_u \to T$.\footnote{To avoid confusion, we change the notation $C$ in Lemma \ref{Lemma 1.1} to $T$.} 
Except for the component $H$, all other components of $\Gamma-\sum_{i=1}^sA'_i-R$ are fiber components of the $\BP^1$-fibration $p_W$
on $W$ induced by $p$. There are two cases depending on whether $H$ is a component of $R$ or not.

\begin{lem}\label{Lemma 3.5}
Suppose that $H$ is not contracted by $\Phi$. Then the following assertions hold.
\begin{enumerate}
\item[{\rm (1)}]
The image $\Phi(H)$ is either $S_1$ or an irreducible component of $D_\infty-S_1$, where $D_\infty$ is the connected component of $D$ such 
that $S_1 \subseteq D_\infty$.
\item[{\rm (2)}]
Suppose that either $G$ is non-cyclic or $G$ is cyclic and $d > 1$. If $\Phi(H)=S_1$ then $\varphi$ preserves the standard 
$\A^1_\ast$-fibration.
\item[{\rm (3)}]
If $G$ is cyclic and $d=1$ then $\varphi$ preserves the standard $\A^1_\ast$-fibration.
\end{enumerate}
\end{lem}
\Proof
(1)~~By Lemma \ref{Lemma 3.4}, $H \not\subset \Supp(R)$ since $H$ is a component of $\Delta$ and is not contracted by 
the morphism $\Phi$. Let $\ell$ be a general fiber of $p$ on $V_u$ and let $C=\Phi(\ell)$. Intersecting $\ell$ with both sides 
of the formula (1), we have $C\cdot(D+K_{V_l})=0$ since $\ell\cdot(\rho^\ast(D+K_{V_u}))=\ell\cdot\left(\Gamma-\sum_iA'_i-R\right)=0$ 
as $H \not\subset \Supp(R)$. Suppose that $G$ is non-cyclic. Then Lemma \ref{Lemma 1.2} shows that 
\begin{eqnarray*}
D+K_{V_l} \sim \ell -(F_1+F_2+F_3) \ge \left(1-\frac{1}{m_1}-\frac{1}{m_2}-\frac{1}{m_3}\right)\ell,
\end{eqnarray*}
where $\ell$ denotes a fiber of the standard $\BP^1$-fibration on $V_l$ by abuse of notation and the difference is an effective 
$\Q$-divisor
\[
L:=\sum_{i=1}^3\frac{1}{m_i}(\ell_i-m_iF_i).
\]
Hence we have
\begin{eqnarray*}
0 \ge \left(1-\frac{1}{m_1}-\frac{1}{m_2}-\frac{1}{m_3}\right)\cdot(C\cdot \ell)
\end{eqnarray*}
because $(C\cdot D+K_{V_l})=0$. Note that $C$ can move as $\dim|\ell|=1$. Since $\left(1-\frac{1}{m_1}-\frac{1}{m_2}-\frac{1}{m_3}\right)<0$, 
it follows that $C\cdot\ell\ge 0$. If $C\cdot\ell=0$, $C$ is a fiber of the fibration $p$, and this is what we have to prove.
 Suppose $C\cdot\ell >0$. Then $L\cdot C > 0$. Hence $C$ meets a boundary component of $\ell_i-m_iF_i$ for some $i=1, 2, 3$. If $C$ meets 
an irreducible component connected to $S_0$, i.e., if $C$ meets the connected component $D_0$ of $D$ containing $S_0$, then the image of $C$ 
by the contraction of $S_0+T_1+T_2+T_3$ is a complete curve contained in $X=\A^2/G$. In fact, $C$ meets $S_0$ because $\Phi$ is the identity 
on $S_0$ and has only two places at infinity (i.e., at the points in the boundary $D$) as so does $\ell$. This is a contradiction. For, 
otherwise, $C$ has at least three places at infinity. Hence $C$ meets an irreducible component of $D_\infty$, and one of the 
components of $D_\infty$ is necessarily $\Phi(H)$. 

(2)~~Suppose that $\Phi(H)=S_1$. Then $C$ has only two places at the boundary over the points $C\cap S_0$ and $C\cap S_1$, both of which 
move by the assumption. This implies that $C$ corresponding to a general $\ell$ does not meet an irreducible component of $L$. 
Hence $C\cdot\ell=0$. 

The case where $G$ is cyclic and $d > 1$ can be handled in a similar way. The $\Q$-divisor $L$ is 
\[
L=\frac{1}{m_0}(\ell_0-m_0F_0)+\frac{1}{m_\infty}(\ell_\infty-m_\infty F_\infty)
\]
in this case.

(3)~~Suppose that $G$ is cyclic and $d=1$. By Lemma \ref{Lemma 1.2}, $D+K_{V_l} \sim -2\ell$. Hence, for $C=\Phi(\ell)$, we have 
\[
0=C\cdot (D+K_{V_l})=-2(C\cdot \ell).
\]
Hence $(C\cdot \ell)=0$. This implies that $C$ is a fiber of the standard $\BP^1$-fibration $p$ on $V_\ell$. Hence $\varphi$ 
preserves the $\A^1_\ast$-fibration.
\QED

\begin{remark}\label{Remark 3.6}{\em
Lemma \ref{Lemma 3.5} shows that the problem becomes harder in the case where $\Phi(H)$ is a component of $L$ or in 
the case where $H$ is contracted by $\Phi$. We do not know if $H$ is not contracted by $\Phi$ provided $H \not\subset \Supp(R)$. 
In fact, Lemma \ref{Lemma 3.4} suggests that if $H$ contracts to an intersection points of two boundary components of $D$ on $V_l$, 
$H$ might not be a component of $R$. Only in the third case of Lemma \ref{Lemma 1.2} where $G$ is cyclic and $d=1$, $H \not\subset R$ 
implies that $H$ is not contracted by $\Phi$. This is because $D$ consists of components $S_0$ and $S_1$.}
\QED
\end{remark}

Let $\ell$ be a general fiber of the standard $\BP^1$-fibration $p : V_l \to C$ and let $B=\Phi^\ast(\ell)$. Since 
$\ell^2=0$ and $\dim |\ell|=1$, it follows that $B^2=0$ and $|B|$ is composed of a linear pencil. Then we have the following result.

\begin{lem}\label{Lemma 3.7}
With the above notations, $B$ is an irreducible curve, and $|B|$ is an irreducible pencil. Further, the following four conditions 
are equivalent. 
\begin{enumerate}
\item[(1)]
$B\cdot(\Delta+K_W)=0$.
\item[(2)]
$B$ is an irreducible smooth rational curve, and $B\cdot R=0$.
\item[(3)]The boundary divisor $\Delta$ contains an irreducible component $\Delta_1$ such that $\Phi$ induces an isomorphism 
between $\Delta_1$ and $S_1$. Further, $B$ does not meet any other irreducible components of $W-\wt{Y}$.
\item[(4)]
$\varphi$ is an automorphism.
\end{enumerate}
\end{lem}
\Proof
Let $Q=S_0\cap\ell$. Since $\wh{\varphi}^{-1}(Q)$ is a single point $\wh{Q}$ on $S_0$ on $\wh{X}_u$ by Lemma \ref{Lemma 2.2} 
and since $\wh{\varphi}$ is \'etale at the point $\wh{Q}$, we have an isomorphism $\wh{\SO}_{W,\wh{Q}} \cong 
\wh{\SO}_{V_l,Q}$. Hence the analytic local behavior of $B$ at $\wh{Q}$ is the same as the one of $\ell$ at $Q$. 
This implies that $B$ is irreducible. 
In fact, since $(B^2)=(\Phi^\ast(\ell))^2=n(\ell^2)=0$ with $n=\deg\Phi$, $|B|$ is composed of an irreducible linear pencil 
$\Lambda$ parametrized by the base curve of $\wt{p} : V_u \to \wt{C}$, where $\wt{C}$ is the normalization of $C$ in 
the function field of $V_u$ (Stein factorization). If $m=\deg_C\wt{C}$, then $B=B_1+\cdots+B_m$ with $B_i \in \Lambda$ and 
hence $S_0\cdot B=\sum_{i=1}^mS_0\cdot B_i\ge m$. Since $\wh{\varphi}|_{S_0}$ is an isomorphism by Lemma \ref{Lemma 2.2}, 
$B$ meets $S_0$ only at $\wh{Q}$. Hence $(S_0\cdot B)_W=(S_0\cdot \ell)_{V_l}=1$, whence we have $m=1$.
\svskip

\noindent
(1)~$\Rightarrow$~(2)~Suppose that $B\cdot(\Delta+K_W)=0$. Since $\Delta$ contains $S_0$ (identified with the proper transform 
on $W$), $B\cdot \Delta >0$. Hence $B\cdot K_W<0$. By the arithmetic genus formula, it follows that $B^2=0$ and $B\cdot K_W=-2$.
Namely $B$ is a smooth rational curve. By the logarithmic canonical divisor formula in Lemma \ref{Lemma 1.2} and the logarithmic ramification 
divisor formula before Lemma \ref{Lemma 3.4}, $B\cdot(\Delta+K_W)= 0$ implies $B\cdot R=0$.
\svskip

\noindent
(2)~$\Rightarrow$~(1)~Since $B\cdot\Phi^\ast(D+K_{V_l})=\deg\Phi\cdot(\ell\cdot D+K_{V_l})=0$ by Lemma \ref{Lemma 1.2}, 
it follows that $B\cdot(\Delta+K_W)=0$ if $B\cdot R=0$.
\svskip

\noindent
(2)~$\Rightarrow$~(3)~Since $B\cdot (\Delta+K_W)=0$, $B\cdot \Delta=2$ and $B\cdot S_0=1$, there exists an irreducible component 
$\Delta_1$ of $\Delta$ such that $B\cdot \Delta_1=1$. Hence $\ell\cdot \Phi_\ast(\Delta_1)=1$. Since $\ell\cdot S_1=1$, 
it follows that $\Phi_\ast(\Delta_1)=S_1$. This implies that $\Phi\vert_{\Delta_1} : \Delta_1 \to S_1$ is an isomorphism. It is 
clear that $B$ does not meet the components of $\Delta-(S_0+\Delta_1)$.
\svskip

\noindent
(3)~$\Rightarrow$~(4)~The restriction $\Phi\vert_B : B \to \ell$ is a finite covering of $\BP^1$ which is totally 
ramifying at $B\cap \Delta_1$ by the condition (3) and \'etale over $B\cap S_0$ by Lemma \ref{Lemma 2.2}. If $\deg\Phi\vert_B
>1$ then $\Phi^{-1}(\ell\cap S_0)$ consists of the points $B\cap S_0$ and the intersection points $\coprod_i(B\cap A_i)$. 
Then one of the components of $\coprod A_i$, say $A_1$, is mapped by $\Phi$ to the curve $S_0$ on $V_l$. In fact, since the point 
$\ell\cap S_0$ moves as $\ell$ moves, $\Phi$ induces a correspondence between $A_1$ and $S_0$ under the correspondence between 
the family $\{B\}$ and $\{\ell\}$. Hence $A_1$ is mapped to the singular point $P_0$ by the normalization morphism $\nu : Y \to X_l$. 
This is a contradiction because $\nu$ is a finite morphism. Thus $\deg\Phi\vert_B=1$. This implies that $\deg \varphi=1$. Namely, 
$\varphi$ is birational. Then the lifting $\wt{\varphi}$ of $\varphi$ onto $\A^2$ is birational. Note that $\wt{\varphi}$ is an automorphism because an \'etale birational endomorphism of $\A^2$ is an automorphism. Hence $\varphi$ is an automorphism. 
\svskip

\noindent
(4)~$\Rightarrow$~(1)~$B$ is a general member of the pencil $|B|$. Hence $B$ is a smooth rational curve isomorphic to $\ell$ by 
the assumption that $\varphi$ is an automorphism. Furthermore, $B$ has two places outside $X^\circ$ as $\ell$ does and $B\cdot S_0=1$. 
Hence $B$ meets one component, say $\Delta_1$, of $\Delta$ in one point transversally and does not meet other components of 
$\Delta-(\Delta_1+S_0)$. Hence $B\cdot(\Delta+K_W)=0$. 
\QED

\section{The case where $G$ is cyclic and $d=1$}
We shall show that GJC holds true under the assumption that $G$ is cyclic and $d=1$. We retain the assumption throughout the present 
section. In view of Theorem \ref{Theorem 2.4} and Lemma \ref{Lemma 3.5}, it remains to show that $H$ may be assumed to be not contracted 
by $\Phi$. 

\begin{lem}\label{Lemma 4.1}
Assume that $\Phi(H)$ is a point $Q_1$ in $V_l$. Then the following assertions hold.
\begin{enumerate}
\item[{\rm (1)}]
$Q_1 \in S_1$.
\item[{\rm (2)}]
Let $C=\Phi(\ell)$ for a general fiber $\ell$ of the standard $\BP^1$-fibration $p : V_u \to T$ \footnote{In the previous sections, 
the base curve is denoted by $C$. In order to avoid confusion, we denote it by $T$ and keep $C$ for the image $\Phi(\ell)$}, 
which is identified with the proper transform on the surface $W$. Then $C \sim aS_1+\ell$ with $a > 0$, where $\ell$ is a fiber of 
the standard $\BP^1$-fibration $p : V_l \to T$.
\item[{\rm (3)}]
Let $\Lambda$ be the linear pencil on $V_u$ consisting of the fibers of $p$. Then $\Phi_\ast(\Lambda)$ is a linear pencil spanned by 
$C$ and $C_0:=aS_1+\ell_0$, where $\ell_0=p^{-1}(p(Q_1))$.
\item[{\rm (4)}]
A general member $C$ has two places at infinity at $Q_0=C\cap S_0$ and $Q_1=C\cap S_1$, where $i(C,S_0;Q_0)=1$, $i(C,S_1;Q_1)=an+1$ 
and $i(C,\ell_0;Q_1)=a$, where $n=|G|$. Hence $C$ has a cusp of type $(an+1,a)$ at $Q_1$.
\item[{\rm (5)}]
Let $\sigma : V' \to V_l$ be a minimal resolution of base points of $\Phi_\ast(\Lambda)$. Let $\sigma'(\Phi_\ast(\Lambda))$ be the 
proper transform of $\Phi_\ast(\Lambda)$ by $\sigma$. Then it is the linear pencil spanned by the proper transform $C':=\sigma'(C)$ 
and the member $C'_0$ obtained from $C_0$, where 
\begin{eqnarray*}
C'_0 &=&  a(S'_1+A_1+\cdots+A_a)+(E_n+E_{n-1}+\cdots+E_1+\ell'_0) \\
     & & +(a-1)B_1+(a-2)B_2+\cdots+B_{a-1},
\end{eqnarray*}
where the last $(-1)$-curve $F$ of $\sigma$, which is a section of $\sigma'(\Phi_\ast(\Lambda))$, meets the curve $B_{a-1}$. 
The weighted dual graph of $\sigma^*(aS_1+\ell_0)$ is given as below.

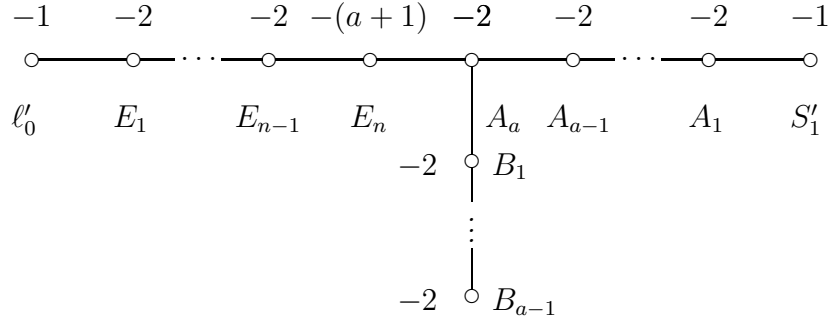
\raisebox{-50mm}{
\begin{picture}(135,60)(25,0)
\setlength{\unitlength}{0.9mm}
\put(10,45){\circle{2}}
\put(11,45){\line(1,0){13}}
\put(25,45){\circle{2}}
\put(26,45){\line(1,0){5}}
\put(32,45){$\dots$}
\put(38,45){\line(1,0){6}}
\put(45,45){\circle{2}}
\put(46,45){\line(1,0){13}}
\put(60,45){\circle{2}}
\put(61,45){\line(1,0){13}}
\put(75,45){\circle{2}}
\put(76,45){\line(1,0){13}}
\put(90,45){\circle{2}}
\put(91,45){\line(1,0){5}}
\put(97,45){$\dots$}
\put(103,45){\line(1,0){6}}
\put(110,45){\circle{2}}
\put(111,45){\line(1,0){13}}
\put(125,45){\circle{2}}
\put(75,44){\line(0,-1){13}}
\put(75,30){\circle{2}}
\put(75,29){\line(0,-1){5}}
\put(75,18){\!$\vdots$}
\put(75,17){\line(0,-1){6}}
\put(75,10){\circle{2}}
\put(7,50){$-1$}
\put(22,50){$-2$}
\put(42,50){$-2$}
\put(72,50){$-2$}
\put(51,50){$-(a+1)$}
\put(72,50){$-2$}
\put(87,50){$-2$}
\put(107,50){$-2$}
\put(122,50){$-1$}

\put(7,35){$\ell'_0$}
\put(22,35){$E_1$}
\put(40,35){$E_{n-1}$}
\put(57,35){$E_n$}
\put(77,35){$A_a$}
\put(86,35){$A_{a-1}$}
\put(107,35){$A_1$}
\put(122,35){$S'_1$}
\put(64,28){$-2$}
\put(64,8){$-2$}
\put(78,28){$B_1$}
\put(78,8){$B_{a-1}$}
\end{picture}}
\item[{\rm (6)}]
The curve $C$ has $n$ cusps of multiplicity $a$ at $Q_1$ and its infinitely near points and no singular points elsewhere. 
Hence the pencil $\sigma'(\Phi_\ast(\Lambda))$ is a linear pencil of $\BP^1$s. 
\item[{\rm (7)}]
Let $\tau$ be the contraction of the curves $S'_1, A_1,\ldots,  A_a, B_1, \ldots,\\ B_{a-1}, E_n, \ldots, E_1$ in this order. 
Then $\tau(V')=\F_n$, $\tau(\ell'_0)$ is a fiber, and $\tau(F)$ is a section disjoint from the minimal section $\tau(S'_0)$. 
The $\BP^1$-fibration on the new completion $\F_n$ of $X_l$ is given by $(\tau\circ\sigma^{-1}\circ\Phi)_*\Lambda$.
\end{enumerate}
\end{lem}

\Proof
(1)~~Othewise, the image $C$ would give a complete curve in the affine surface $X=\A^2/G$. 

(2)~~Since $V_l \cong \F_n$, $C$ is written as $C \sim aS_1+b\ell$. Since $C\cdot S_0=1$, it follows that $b=1$. 

(3)~~As $\ell$ moves in the pencil, $C$ is parametrized by the point $C\cap S_0$. Hence $\Phi_*(\Lambda)$ is a linear pencil 
parametrized by $S_0$. 

(4)~~$C$ has only one place over each of $S_0\cap C$ and $S_1\cap C$. So, the assertion follows straightforwardly.

(5)~~After the $n$ blowing-ups with centers at the intersection points of the proper transforms of $S_1$ and $C$, we have 
a linear equivalence relation
\[
aS_1^{(n)}+\ell_0^{(n)}+E_1+\cdots+E_n \sim C^{(n)},
\]
where $S_1^{(n)}, \ell_0^{(n)}, C^{(n)}$ are the proper transforms of $S_1, \ell_0, C$ respectively. Now $C^{(n)}$ is a smooth 
curve with $(C^{(n)})^2=(C^2)-a^2n=2a$, $C^{(n)}$ touches $E_n$ in one point, say $Q_n$, with multiplicity $i(C^{(n)},E_n;Q_n)
=a$. We blow up $Q_n$ and its infinitely near points $a$ times to separate $C^{(n)}$ from $E_n$. This process gives the 
exceptional curves $A_1, \ldots, A_a$. Since $(S_1^{(n)})^2=n-n=0$, we have the upper horizontal line of the dual graph in the 
statement. Further we have a linear equivalence relation 
\[
a(S_1^{(n+a)}+A_1+\cdots+A_a)+\ell_0^{(n+a)}+E_1+\cdots +E_n \sim C^{(n+a)}.
\]
Since the intersection point $A_a\cap C^{(n+a)}$ is a base point of the linear pencil $\Lambda^{(n+a)}$ which is the proper 
transform of $\Lambda$, we need $a$ times more blowing-ups with centers at $A_a\cap C^{(n+a)}$ and its infinitely near points, 
which gives a vertical chain $B_1+\cdots+B_{a-1}$ in the graph. The last exceptional curve $B_a (=F$ in the statement) is 
a section of the linear pencil $\sigma'(\Phi_*(\Lambda))$. 

The verification of the assertions (6) and (7) is straightforward. 
\QED

Replace the original surface $V_l$ by this birational transform in the assertion (7) as a new completion of $\wh{X}_l$ 
such that the pullback of the standard $\BP^1$-fibration has no base points on $V_u$. Then $\Phi$ is unchanged in the interior 
part $\wh{X}_u$ of $W$. Hence we may assume that $\Phi(H)$ is not contracted. Theorem \ref{Theorem 2.4} and Lemma \ref{Lemma 3.5} 
show that $\varphi$ is an automorphism in the case $G$ is cyclic and $d=1$. 
\svskip

\noindent
{\bf Proof of Theorem in Introduction.}~~
As explained in Lemma \ref{Lemma 2.1} (3), some power $\varphi^N$ lifts up to a $G$-equivariant \'etale endomorphism 
$\wt{\varphi}^N$ of $\A^2$, and $\varphi$ is an automorphism if and only if $\varphi^N$ is an automorphism of $X$. Further,
 $\varphi^N$ is an automorphism if and only if $\wt{\varphi}^N$ is an automorphism of $\A^2$. So, we may assume that $\varphi$ 
lifts up to a $G$-equivariant \'etale endomorphism $\wt{\varphi}$ of $\A^2$. By the hypothesis, $G$ has a subgroup $H$ of order 
$2$. Since $G$ is small, so is $H$. Hence $H$ acts on $\A^2$ as $(x,y) \mapsto (-x,-y)$. Moreover, $\wt{\varphi}$ induces 
an \'etale endomorphism $\varphi' : Y \to Y$, where $Y=\A^2/H$. By the last result which we have proved for a cyclic group $G$ 
and $d=1$, $\varphi'$ and hence $\wt{\varphi}$ are automorphisms.
\QED


\begin{thebibliography}{15}

\bibitem{Artin}
M. Artin, {\em On isolated rational singularities of surfaces}, Amer. J. Math. {\bf 88} (1966), 129--136.

\bibitem{Ax}
J. Ax, {\em Injective endomorphisms of varieties and schemes}, Pacific J. Math. {\bf 31} (1969), 1--7.

\bibitem{B}
E. Brieskorn, {\em Rationale Singularit\"aten komplexer Fl\"achen}, Invent. Math. {\bf 4} (1968), 336--358.

\bibitem{Chevalley}
C. Chevalley, {\em Invariants of finite groups generated by reflections}, Amer. J. Math. {\bf 77} (1955), 778--782.

\bibitem{DP}
A. Dubouloz and K. Palka, {\em The Jacobian conjecture fails for pseudo-planes}, Adv. Math. {\bf 339} (2018), 248--284.

\bibitem{EGA}
A. Grothendieck, {\em \'El\'ements de g\'eom\'etrie alg\'ebrique, IV. \'Etude locale des sch\'emas et des morphismes de sch\'emas, III}, 
Inst. Hautes \'Etudes Sci. Publ. Math. No. {\bf 28} (1966), 255 pp.

\bibitem{GM}
R.V. Gurjar and M. Miyanishi, {\em Jacobian problem for singular surfaces}, J. Math. Kyoto Univ. {\bf 48} (2008), 
757--764.

\bibitem{Iitaka}
S. Iitaka, {\em Algebraic Geometry. An introduction to Birational Geometry of Algebraic Varieties}, Graduate Texts in Mathematics, 
{\bf 76}, Springer-Verlag, New York-Berlin, 1982. x+357 pp. 

\bibitem{Lipman}
J. Lipman, {\em Rational singularities, with applications to algebraic surfaces and unique factorization}, 
Inst. Hautes Études Sci. Publ. Math. No. {\bf 36} (1969), 195--279.

\bibitem{MM}
K. Masuda anad M. Miyanishi, {\em \'Etale endomorphisms of algebraic surfaces with $G_m$-actoions}, Math. Ann. {\bf 319} (2001),
493--516.

\bibitem{M1}
M. Miyanishi, {\em \'Etale endomorphisms of algebraic varieties}, Osaka J. Math. {\bf 22} (1985), 345--364.

\bibitem{M2}
M. Miyanishi, {\em Non-complete Algebraic Surfaces}, Lecture Notes in Math. {\bf 857}, Springer Verlag, 1981.

\bibitem{M3}
M. Miyanishi, {\em Open Algebraic Surfaces}, CRM Monograph Series {\bf 12}, Amer. Math. Soc. , 2001.

\bibitem{KM}
M. Miyanishi and K. Masuda, {\em Generalized Jacobian conjecture and related topics}, Proc. Internat. Colloq. on 
Algebra, Arithmetic and Geometry, 2002, 427--466, Tata Institute for Fundamental Research and International 
Mathematical Union.

\bibitem{KM2}
M. Miyanishi and K. Masuda, {\em Affine pseudo-planes with torus actions}, Transform. Groups {\bf 11} (2006), no. 2, 249--267.

\bibitem{Tyurina}
G.N. Tyurina, {\em Absolute isolatedness of rational singularities and triple rational points}, Functional Analysis and 
its applications {\bf 2}(1968), no.4, 324-333. 

\end{thebibliography}
\end{document}